\documentclass[letterpaper, 10 pt, conference]{ieeeconf}  

\IEEEoverridecommandlockouts                              

\usepackage{times}
\usepackage{cite}
\usepackage[hidelinks]{hyperref}
\usepackage{enumerate}
\usepackage{mathtools}
\usepackage{dsfont}
\usepackage{amsmath}
\usepackage{amssymb} 
\usepackage{color}
\usepackage{gensymb}
\usepackage{xcolor}
\usepackage{multirow}
\newtheorem{theorem}{Theorem}

\newtheorem{problem}{Problem}

\newtheorem{remark}{Remark}
\newtheorem{definition}{Definition}
\newtheorem{assumption}{Assumption}
\newtheorem{ex}{Example}

\title{\LARGE \bf
Risk-Minimizing Two-Player Zero-Sum Stochastic\\ Differential Game via Path Integral Control}

\author{Apurva Patil$^1$ \and Yujing Zhou$^2$ \and David Fridovich-Keil$^3$ \and Takashi Tanaka$^3$ 
\thanks{$^{1}$Walker Department of Mechanical Engineering, University of Texas at Austin, {\tt\small apurvapatil@utexas.edu}.
$^{2}$Department of Mechanical and Aerospace Engineering, Princeton University, {\tt\small yz1324@princeton.edu}.
$^{3}$Department of Aerospace Engineering and Engineering Mechanics, University of Texas at Austin, {\tt\small dfk@utexas.edu}, {\tt\small ttanaka@utexas.edu}.}
}
\begin{document}

\maketitle
\thispagestyle{empty}
\pagestyle{empty}

\begin{abstract}
This paper addresses a continuous-time risk-minimizing two-player zero-sum stochastic differential game (SDG), in which each player aims to minimize its probability of failure. Failure occurs in the event when the state of the game enters into predefined undesirable domains, and one player's failure is the other's success. We derive a sufficient condition for this game to have a saddle-point equilibrium and show that it can be solved via a Hamilton-Jacobi-Isaacs (HJI) partial differential equation (PDE) with Dirichlet boundary condition. Under certain assumptions on the system dynamics and cost function, we establish the existence and uniqueness of the saddle-point of the game. We provide explicit expressions for the saddle-point policies which can be numerically evaluated using path integral control. This allows us to solve the game online via Monte Carlo sampling of system trajectories. We implement our control synthesis framework on two classes of risk-minimizing zero-sum SDGs: a disturbance attenuation problem and a pursuit-evasion game. Simulation studies are presented to validate the proposed control synthesis framework.
\end{abstract}
\section{Introduction}\label{Sec: Introduction}
Interactions among multiple agents are prevalent in many fields such as economics, politics, and engineering. Game theory studies the collective decision-making process of multiple interacting agents \cite{bacsar1998dynamic}. Two-person zero-sum games involve two players with conflicting interests, and one player's gain is the other's loss. Pursuit-evasion games (competition between a pursuer and an evader) \cite{nahin2012chases,7172219} and robust control  (competition between a controller and the nature) \cite{bacsar1998dynamic} are some examples of two-player zero-sum games. In this paper, we consider a two-player stochastic differential game (SDG) in which the outcome of the game depends not only on the decision of both players but also on the stochastic input added by the nature.\par

When the game dynamics and cost functions are known, the saddle-point equilibrium of a two-player zero-sum SDG can be characterized by the Hamilton-Jacobi-Isaacs (HJI) partial differential equation (PDE). The analytical solutions of the HJI PDEs are in general not available, and one needs to resort on numerical methods such as grid-based approaches \cite{falcone2006numerical}, \cite{huang2014automation} to solve these PDEs approximately. 
However, the grid-based approaches suffer from \textit{curse of dimensionality}, making them computationally intractable for systems with large dimensions \cite{mitchell2005time}. Moreover, in general, the solutions can not be computed in real-time using these methods; they need to be precomputed as lookup tables and recalled for the use in an online setting \cite{huang2014automation}.
Several reinforcement learning algorithms have also been proposed to find approximate solutions to game problems. A reinforcement-learning-based adaptive dynamic programming algorithm is proposed in \cite{vrabie2011adaptive} to determine online a saddle-point solution of linear continuous-time two-player zero-sum differential games. A deep reinforcement learning algorithm based on updating players' policies simultaneously is proposed in \cite{prajapat2021competitive} to solve two-player zero-sum games.
These methods assume deterministic game dynamics and do not consider system uncertainties. In the presence of system uncertainties, the performance and safety of both players are affected unpredictably, if the uncertainties are not accommodated while designing policies. An effective uncertainty evaluation method, the multivariate probabilistic collocation, was used in \cite{liu2020adaptive} with integral reinforcement learning to solve multi-player SDGs for linear system dynamics online. A two-person zero-sum stochastic game with discrete states and actions is solved in \cite{lin2017multiagent} using Bayesian inverse reinforcement learning. Common challenges in the learning-based methods include training efficiency, rigorous theoretical guarantees on convergence and optimality. Moreover, these approaches do not explicitly take into account the players' failure probabilities while synthesizing their policies. \par

In our work, we formulate a continuous-time, nonlinear, two-player zero-sum SDG on a state space modeled by an It\^o stochastic differential equation. Since the stochastic uncertainties in our model are unbounded, both players have nonzero probabilities of failure. Failure occurs when the state of the game enters into predefined undesirable domains, and one player's failure is the other's success. Our objective is to solve a game in which each player seeks to minimize its risk of failure (failure probability)\footnote{Throughout this paper the word ``risk" simply means the probability of failure. It is not our intention to discuss various risk measures existing in the literature (e.g. \cite{artzner1999coherent, dixit2022risk}).} along with its control cost; hence, the name \textit{risk-minimizing zero-sum SDG}. We explain the risk-minimizing zero-sum SDG via the following example:
\begin{ex}
Consider a pursuit-evasion game in which the pursuer catches the evader if they are less than a certain distance $\rho$ away from each other. In this setting, the evader wishes to minimize its probability of entering the ball of radius $\rho$ centered at the pursuer's location. Whereas, the pursuer wishes to minimize the probability of staying out of the ball of radius $\rho$ centered at the evader's location. The goal of each player is to balance the trade-off between the above probabilities (probabilities of failure) and the control cost (for e.g., their energy consumption). 
This problem can be formulated as a risk-minimizing zero-sum SDG. 
\end{ex}\par
We derive a sufficient condition for this game to have a saddle-point equilibrium, and show that it can be solved via an HJI PDE with \textit{Dirichlet boundary condition}. Under certain assumptions on the system dynamics and cost function, we establish existence and uniqueness of the saddle-point equilibrium of the formulated risk-minimizing zero-sum SDG. Furthermore, explicit expressions for the saddle-point policies are derived which can be numerically evaluated using path integral control. The idea behind the path integral control is to use the Feynman-Kac lemma \cite{oksendal2013stochastic} and solve a linear PDE via Monte Carlo samples of system trajectories \cite{kappen2005path}, \cite{williams2017model}. The Monte Carlo simulations can be massively parallelized through the use of graphics processing units (GPUs); hence this approach is less susceptible to curse of dimensionality. The use of path integral technique to solve stochastic games was proposed in \cite{vrushabh2020robust}. In this paper, we generalize their work and develop a path integral formulation to solve HJI PDEs with Dirichlet boundary conditions and find saddle-point equilibria of risk-minimizing zero-sum SDGs. The proposed framework allows us to solve the game online using Monte Carlo simulations of system trajectories, without the need of any offline training or precomputations.  \par
The contributions of this work are as follows: 1) We formulate a continuous-time risk-minimizing zero-sum SDG in which players aim at balancing the trade-off between the failure probability and control cost. A sufficient condition for this game to have a saddle-point equilibrium is derived, and it is shown that this game can be solved via an HJI PDE with Dirichlet boundary condition.
2) Under certain assumptions on the system dynamics and cost function, we establish the existence and uniqueness of the saddle-point solution. We also obtain explicit expressions for the saddle-point policies which can be numerically evaluated using path integral control. 3) The proposed control synthesis framework is validated by applying it on two classes of risk-minimizing zero-sum SDGs, namely a disturbance attenuation problem and a pursuit-evasion game.  

\subsection*{Notation}
Bold symbols such as $\pmb{x}$ represent random variables. If a stochastic process $\pmb{x}(s), s\geq t$ starts from $x$ at time $t$, then let $P_{x,t}\left(\mathcal{E}\right)$ denote the probability of event $\mathcal{E}$ conditioned on $\pmb{x}(t)\!=\!x$, and let $\mathbb{E}_{x, t}\left[F\left(\pmb{x}\right)\right]$ denote the expectation of a functional $F\left(\pmb{x}\right)$ conditioned on $\pmb{x}\left(t\right)\!=\!x$. Let $\mathds{1}_{\mathcal{E}}$ be an indicator function, that returns $1$ when the condition $\mathcal{E}$ holds and 0 otherwise. $\text{Tr}(A)$ denotes the trace of a matrix $A$.

\section{Problem Formulation}\label{Sec: Problem Formulation}
We consider a two-player zero-sum stochastic differential game (SDG) on a finite time horizon $t\in[t_0, T]$, $t_0<T$. Consider a class of control-affine stochastic systems described by the following It\^{o} stochastic differential equation (SDE):

\begin{equation}\label{SDE}
\begin{aligned}
  d\pmb{x}(t)=&{f}\left(\pmb{x}(t), t\right)dt+{G_u}\left(\pmb{x}(t), t\right){u}(\pmb{x}(t), t)dt\\
  &\!\!\!+ G_v\left(\pmb{x}(t), t\right) v\left(\pmb{x}(t), t\right)dt+{\Sigma}\left(\pmb{x}(t), t\right)d\pmb{w}(t)
\end{aligned}
\end{equation}
where $\pmb{x}(t)\in\mathbb{R}^n$ is the state, ${u}(\pmb{x}(t), t)\in\mathbb{R}^m$ is the control input of the first player (henceforth called the agent), and ${v}(\pmb{x}(t), t)\in\mathbb{R}^l$ that of the second player (called the adversary). $\pmb{w}(t)\in\mathbb{R}^k$ is a $k$-dimensional standard Wiener process on a suitable probability space $\left(\Omega, \mathcal{F}, P\right)$. We assume sufficient regularity in the functions ${f}\left(\pmb{x}(t), t\right)\in\mathbb{R}^n$, ${G_u}\left(\pmb{x}(t), t\right)\in \mathbb{R}^{n\times m}$, ${G_v}\left(\pmb{x}(t), t\right)\in \mathbb{R}^{n\times l}$ and ${\Sigma}\left(\pmb{x}(t), t\right)\in\mathbb{R}^{n\times k}$ so that a unique strong solution  of (\ref{SDE}) exists \cite{oksendal2013stochastic}. 
Both the control inputs $u$, and $v$ are assumed to be square integrable (i.e., of finite energy). In the rest of the paper, for notational compactness, the functional dependencies on $x$ and $t$ are dropped whenever it is unambiguous. \par

Let $\mathcal{X}_{s}\subseteq\mathbb{R}^n$ be a bounded open set representing a safe region, $\partial\mathcal{X}_{s}$ be its boundary, and closure $\overline{\mathcal{X}_{s}}=\mathcal{X}_{s}\cup\partial\mathcal{X}_{s}$. Suppose that the agent tries to keep the system (\ref{SDE}) in the safe set $\mathcal{X}_s$ for the entire time horizon $[t_0, T]$ of the game, whereas the adversary seeks the opposite. For example, in pursuit-evasion games, the safe set $\mathcal{X}_s$ could be a region outside the ball of radius $\rho$, centered at the adversary's location. Or, in the disturbance rejection problems, if an agent wishes to navigate through obstacles in the presence of adversarial disturbances, then the region outside obstacles could be considered as a safe set. Suppose, when the game starts at $t_0$, the system is in the safe set i.e., $\pmb{x}(t_0)=x_0\in{\mathcal{X}_s}$. If the system leaves the region $\mathcal{X}_{s}$ at any time $t\in(t_0, T]$, we say that the agent fails. On the other hand, the adversary fails if the system stays in $\mathcal{X}_{s}$ for all $t\in[t_0, T]$. Therefore, we define the agent's probability of failure $P^{\mathrm{ag}}_{\mathrm{fail}}$ as   
\begin{equation}\label{pfail}
    P^{\mathrm{ag}}_{\mathrm{fail}}\!\coloneqq\!P_{x_0,t_0}\!\!\left(\bigvee_{t\in(t_0, T]} \!\!\!\!\pmb{x}(t)\notin \mathcal{X}_{s}\!\!\right)
\end{equation}
and the adversary's probability of failure $P^{\mathrm{ad}}_{\mathrm{fail}} \coloneqq 1-P^{\mathrm{ag}}_{\mathrm{fail}}$. The $\bigvee$ symbol represents a logical OR implying existence of a satisfying event among a collection. We define the terminal time $\pmb{t}_{f}$ of the game as
\begin{equation}\label{tf}
\pmb{t}_{f} \coloneqq 
\begin{cases}
T, & \!\!\!\!\!\!\!\!\!\!\!\!\!\!\!\!\!\!\!\!\!\!\!\!\!\!\!\!\!\!\!\!\!\!\!\!\!\!\text{if}\;\; \pmb{x}(t)\in\mathcal{X}_{s}, \forall t\in(t_0, T),\\
\text{inf}\;\{t\in(t_0, T) : \pmb{x}(t)\notin\mathcal{X}_{s}\}, & \text{otherwise}.
 \end{cases}
\end{equation}
Alternatively, $\pmb{t}_f$ can be defined as 
\begin{equation}\label{tf with Q}
 \pmb{t}_{f} \coloneqq \text{inf}\{t> t_0: (\pmb{x}(t), t)\notin \mathcal{Q}\} 
\end{equation}
where $\mathcal{Q}=\mathcal{X}_s\times[t_0, T)$ is a bounded set with the boundary  $\partial\mathcal{Q}=\left(\partial\mathcal{X}_s\times[t_0,T]\right)\cup\left(\mathcal{X}_s\times\{T\}\right)$, and closure $\overline{\mathcal{Q}}=\mathcal{Q}\cup\partial\mathcal{Q}=\overline{\mathcal{X}_s}\times[t_0, T]$. Note that by the above definitions, $\left(\pmb{x}(\pmb{t}_{f}),\pmb{t}_{f}\right)\in\partial{\mathcal{Q}}$ and the agent's failure probability $P^{\mathrm{ag}}_{\mathrm{fail}}$ in (\ref{pfail}) can be written in terms of $\pmb{t}_f$ as
\begin{equation}\label{pfail2}
    \!P_{x_0,t_0}\!\!\left(\bigvee_{t\in(t_0, T]} \!\!\!\!\pmb{x}(t)\notin \mathcal{X}_{s}\!\!\right)\!=\!\mathbb{E}_{x_0, t_0}\left[\mathds{1} _{\pmb{x}(\pmb{t}_{f})\in \partial\mathcal{X}_{s}}\right]\!.
\end{equation}
Since the stochastic uncertainty of system (\ref{SDE}) is modeled with an unbounded distribution, both the agent and the adversary have nonzero probabilities of failure. In our two-player SDG setting, we assume that both players aim to design an optimal policy against the \textit{worst} possible opponent's policy such that their own \textit{risk of failure} is minimized. Therefore, we define the following \textit{risk-minimizing} cost function:
\begin{equation}\label{C}
\begin{split}
 C\left(x_0, t_0; {u}, {v}\right)&\coloneqq\eta\,\mathbb{E}_{x_0, t_0}\left[\mathds{1} _{\pmb{x}(\pmb{t}_{f})\in \partial\mathcal{X}_{s}}\right]\\
 &\!\!\!\!\!\!\!\!\!\!\!\!\!\!\!\!\!\!\!\!\!\!\!\!\!\!\!\!\!\!\!\!\!\!\!\!\!\!\!\!\!\!\!+\mathbb{E}_{x_0, t_0}\!\Bigg[\!\psi\!\left(\pmb{x}(\pmb{t}_{f}\!)\right)\!\cdot\!\mathds{1} _{\pmb{x}(\pmb{t}_{f})\in \mathcal{X}_{s}}\!+\!\!\int_{t_0}^{\pmb{t}_{f}}\!\!\!\!L\!\left(\pmb{x}(t), \pmb{u}(t), \pmb{v}(t), t\right)\!dt\!\Bigg]\!.\\
\end{split}
\end{equation}
The first term indicates the penalty associated with the agent's failure with the weight parameter $\eta>0$. $\psi\left(\pmb{x}(\pmb{t}_{f})\right)$ and $L\left(\pmb{x}(t), \pmb{u}(t), \pmb{v}(t), t\right)$ denote the terminal and running costs, respectively. Note that the game ends at $\pmb{t}_f$ and the system doesn't evolve after that (this is motivated by the applications where ``collision" or ``capture" ends the game). Therefore, the running cost is integrated over the time horizon $[t_0, \pmb{t}_f]$. The agent tries to minimize $C$ by controlling $u$, whereas the adversary tries to maximize it by controlling $v$. The weight parameter $\eta$ balances the trade-off between the control cost and the failure probability. \par
Notice that if we define $\phi:\overline{\mathcal{X}_s}\to\mathbb{R}$ as\footnote{In the sequel, function $\phi(x)$ sets a boundary condition for a PDE. In order to guarantee the existence of a solution of such a PDE, assumptions on the regularity of $\phi(x)$ (e.g., continuity on $\overline{\mathcal{X}_s}$) are often required. When these requirements are necessary, (\ref{phi(x)}) can be approximated using a smooth bump function $B(x)$ as $\phi(x) \approx \psi(x)B(x) + \eta\left(1-B(x)\right)$.}: 
\begin{equation}\label{phi(x)}
    \phi\left({x}\right)\coloneqq\psi\left({x}\right)\cdot \mathds{1} _{{x}\in \mathcal{X}_{s}}+\eta\cdot\mathds{1} _{{x}\in \partial\mathcal{X}_{s}},
\end{equation}
then, the first term in (\ref{C}) can be absorbed in a new terminal cost function $\phi$ as follows:
\begin{equation}\label{C with phi}
\begin{aligned}
    {C}\!\left(x_0,t_0;{u}, {v}\right)\!=\!\mathbb{E}_{x_0, t_0}\!\Big[\phi\left(\pmb{x}(\pmb{t}_{f})\right)\!+\!\!\int_{t_0}^{\pmb{t}_{f}}\!\!\!\!L\!\left(\pmb{x}, \pmb{u}, \pmb{v}, t\right)dt\Big]\!.
\end{aligned}
\end{equation}
In this paper, we consider the following running cost that is quadratic in $u$ and $v$:
\begin{equation}\label{running cost}
\begin{aligned}
 L\!\left(\pmb{x}, \pmb{u}, \pmb{v}, t\right)\!=\!V\!\left(\pmb{x}, t\right)\!+\!\frac{1}{2}\pmb{u}^{\!T}\!{R_u}\!\left(\pmb{x}, t\right)\pmb{u}\!-\!\frac{1}{2}\pmb{v}^{\!T}\!{R_v}\!\left(\pmb{x}, t\right)\pmb{v}\!
\end{aligned}
\end{equation}
where $V\left(\pmb{x}, t\right)$ denotes a state dependent cost, and $R_u\left(\pmb{x}, t\right)\in\mathbb{R}^{m\times m}$ and $R_v\left(\pmb{x}, t\right)\in\mathbb{R}^{l\times l}$ are given positive definite matrices (for all values of $\pmb{x}$ and $t$). Now, we formulate our risk-minimizing zero-sum SDG as follows: 

\begin{problem} [Risk-Minimizing Zero-Sum SDG]\label{Problem: Risk-minimizing SOC problem}

\begin{align}\label{risk-minimizing SOC problem (tf)}
\min_{u} \max_{v} \; & \mathbb{E}_{x_0, t_0}\!\!\left[\!\phi\left(\pmb{x}(\pmb{t}_{f})\right)\!+\!\!\!\int_{t_0}^{\pmb{t}_{f}}\!\!\!\left(\!\frac{1}{2}\pmb{u}^\top\!\!R_u\pmb{u}\!-\!\frac{1}{2}\pmb{v}^\top\!\!R_v\pmb{v}\!+\!V\!\!\right)\!dt\!\right] \notag\\
\textrm{s.t.} \;\; d\pmb{x} = & {f} dt +{G_u}{u}dt+ G_vvdt + {\Sigma}d\pmb{w},\\
&\pmb{x}(t_0)=x_0. \notag
\end{align}   
where the admissible policies $u$, $v$ are measurable with respect to the $\sigma$-algebra generated by $\pmb{x}(s), t_0\leq s\leq t$.
\end{problem}
Note that Problem \ref{Problem: Risk-minimizing SOC problem} is a \textit{variable-terminal-time} zero-sum SDG where the terminal time is determined by (\ref{tf with Q}).

\section{Synthesis of Minimax Policies}\label{Section: Main Results}
This section presents the main results of the paper. In Section \ref{Section: Risk-Minimizing HJI PDE}, we show that Problem \ref{Problem: Risk-minimizing SOC problem} can be solved via an HJI PDE with appropriate Dirichlet boundary condition. In Section \ref{Sec: PI}, we find a solution of a class of risk-minimizing zero-sum SDGs via path integral control. 
\subsection{HJI PDE with Dirichlet Boundary Condition}\label{Section: Risk-Minimizing HJI PDE}
Notice that the cost function of the risk-minimizing zero-sum SDG (\ref{risk-minimizing SOC problem (tf)}) possesses the time-additive Bellman structure. Therefore, Problem \ref{Problem: Risk-minimizing SOC problem} can be solved by utilizing the principle of dynamic programming. For each $(x,t)\in\overline{\mathcal{Q}}$, and admissible policies $u$, $v$ over $[t, T)$, define the cost-to-go function:
\begin{equation}\label{value function}
\begin{aligned}
 {C}\left(x, t; {u},{v}\right)=&\mathbb{E}_{x, t}\big[\phi\left(\pmb{x}(\pmb{t}_{f})\right)\big]\\
 &\!\!\!\!\!\!\!\!\!\!\!\!\!\!\!\!\!\!\!\!\!\!\!\!+\mathbb{E}_{x, t}\left[\int_{t}^{\pmb{t}_{f}}\left(\frac{1}{2}\pmb{u}^\top R_u\pmb{u}-\frac{1}{2}\pmb{v}^\top R_v\pmb{v}+V\right)\!dt\right].  
\end{aligned}
\end{equation}

\begin{definition}[Saddle-point solution]\label{Def: SP}\cite[Chapter 2]{bacsar2008h}:
Given a two-player zero-sum differential game, a pair of admissible policies $(u^*, v^*)$ over $[t, T)$ constitutes a saddle-point solution, if for each $(x,t)\in\overline{\mathcal{Q}}$, and admissible policies $(u, v)$ over $[t, T)$,
\begin{equation*}
    C(x,t;u^*,v)\leq C^* \coloneqq C(x,t;u^*,v^*) \leq C(x,t;u,v^*).
\end{equation*}
The quantity $C^*$ is the \textit{value} of the game. The value of the game is defined if it satisfies the following relation
\begin{equation*}\label{Isaacs condition}
\begin{aligned}
 C^* \!\!=\! \min_{u} \max_{v} {C}\left(x, t; {u}, {v}\right) = \max_{v} \min_{u} {C}\left(x, t; {u}, {v}\right).  
\end{aligned}
\end{equation*}
\end{definition}
\vspace{2mm}
The following theorem provides the sufficient condition for a saddle-point solution of Problem \ref{Problem: Risk-minimizing SOC problem} to exist.
\begin{theorem}\label{theorem: solution to risk-minimizing soc}
Suppose there exists a function $J:\overline{\mathcal{Q}}\rightarrow \mathbb{R}$ such that 
\begin{enumerate}[(a)]
    \item $J(x,t)$ is continuously differentiable in $t$ and twice continuously differentiable in $x$ in the domain $\mathcal{Q}$;
    \item $J(x,t)$ solves the following stochastic HJI PDE:
    \begin{equation}\label{HJB PDE}
  \!\!\!\!\!\!\!\!\!\!\begin{cases}
     \begin{aligned}
         \!\!-\partial_tJ\!=&V\! +\!f^\top\!\partial_xJ\!+\!\frac{1}{2}\text{Tr}\left(\Sigma\Sigma^\top\partial^2_xJ\right)\\ &\!\!\!\!\!\!\!\!\!\!\!\!\!\!\!\!\!\!\!+\frac{1}{2}\!\left(\partial_xJ\right)^\top\!\!\left(G_vR_v^{-1}G_v^\top\!-\!G_uR_u^{-1}G_u^\top\right)\!\partial_xJ,
          \end{aligned} &  \forall(x,t)\!\in\!\mathcal{Q}, \\
          \vspace{-4mm}&\\
    \!\!\underset{\substack{(x,t)\to(y,s) \\ (x,t)\in\mathcal{Q}}}{\lim}J(x,t)=\phi(y), & \!\!\!\forall(y,s)\in\partial\mathcal{Q}.
  \end{cases}
    \end{equation}
\end{enumerate}
Then, the following statements hold:
\begin{enumerate}[(i)]
\item $J(x,t)$ is the value of the game formulated in Problem \ref{Problem: Risk-minimizing SOC problem}. That is,
\begin{equation}\label{J as value function}
   \begin{aligned}
     J\left(x, t\right)= &\min_{u} \max_{v} {C}\left(x, t; {u}, {v}\right) \\
        = &\max_{v} \min_{u} {C}\left(x, t; {u}, {v}\right),\;\:\forall\;(x,t)\!\in\!\overline{\mathcal{Q}}.
   \end{aligned}
\end{equation}
\item The optimal solution to Problem \ref{Problem: Risk-minimizing SOC problem} is given by 
\begin{equation}\label{u*}
\begin{aligned}
  u^*(x,t)&=-R_u^{-1}\!\left(x, t\right){G_u}^\top\!\!\left(x, t\right)\partial_xJ\!\left(x, t\right), 
\end{aligned}
\end{equation}
\begin{equation}\label{v*}
     v^*(x,t)=R_v^{-1}\!\left(x, t\right){G_v}^\top\!\!\left(x, t\right)\partial_xJ\!\left(x, t\right).
\end{equation}
\end{enumerate}
\end{theorem}
\vspace{2mm}
\begin{proof}
See Appendix A.
\end{proof}

\begin{remark}
Theorem \ref{theorem: solution to risk-minimizing soc} does not say anything about the existence of a function $J(x,t)$ satisfying statements (a) and (b), and it is not in the scope of this paper. However, in Section \ref{Sec: PI}, we focus on a special case in which (\ref{HJB PDE}) can be linearized where the existence and uniqueness of such a function is guaranteed. 
\end{remark}\par

\subsection{Path Integral Formulation}\label{Sec: PI}
 In this section, we derive a path integral formulation to solve a class of risk-minimizing zero-sum SDGs that satify certain assumptions on the system dynamics and cost function. Let $\xi(x,t)$ be the logarithmic transformation (known as Cole-Hopf transformation in the PDE literature) of the value function $J(x,t)$ defined as
\begin{equation}\label{exp transformation}
 J(x,t) = -\lambda\,\text{log}\left(\xi\left(x,t\right)\right)
\end{equation}
where $\lambda$ is a proportionality constant to be defined. Applying the transformation in (\ref{exp transformation}) to (\ref{HJB PDE}) yields

 \begin{equation}\label{transformed HJB PDE}
  \!\!\!\!\!\!\begin{cases}
     \begin{aligned}
         \!\partial_t\xi\!=&\frac{V\xi}{\lambda}-\!\frac{1}{2}\text{Tr}\!\left(\Sigma\Sigma^\top\!\partial^2_x\xi\right)\!+\!\frac{1}{2\xi}\!\left(\partial_x\xi\right)^{\!T}\!\!\Sigma\Sigma^\top\!\partial_x\xi\\
         &\!\!\!\!\!\!\!\!\!\!\!\!\!\!\!\!+\!\frac{\lambda}{2\xi}\!\left(\partial_x\xi\right)^\top\!\!\left(G_vR_v^{-1}G_v^\top\!-\! G_uR_u^{-1}G_u^\top\right)\!\partial_x\xi\!-\!f^\top\partial_x\xi,
          \end{aligned} &  \\
          & \!\!\!\!\!\!\!\!\!\!\!\!\!\!\!\!\!\!\!\!\!\!\!\!\!\!\!\!\!\!\forall(x,t)\in\mathcal{Q}, \\
    \!\!\underset{\substack{(x,t)\to(y,s) \\ (x,t)\in\mathcal{Q}}}{\lim}\xi(x,t)\!=\!\text{exp}{\left(-\frac{\phi(y)}{\lambda}\right)}, & \!\!\!\!\!\!\!\!\!\!\!\!\!\!\!\!\!\!\!\!\!\!\!\!\!\!\!\!\!\!\!\!\!\!\forall(y,s)\in\partial\mathcal{Q}.
  \end{cases}
    \end{equation}
Now, we make the following assumption:
\begin{assumption}\label{Assumption: linearity}
For all $(x,t)\in\overline{\mathcal{Q}}$, there exists a constant $\lambda>0$ such that 
\begin{equation}\label{lambda}
 \begin{aligned}
    \Sigma(x, t)\Sigma^\top(x, t) \!= &\lambda G_u(x, t)R_u^{-1}(x,t) G_u^\top(x, t)\\
    &-\lambda G_v(x, t)R_v^{-1}(x,t) G_v^\top(x, t).\\
 \end{aligned}
\end{equation}
\end{assumption}
\vspace{1mm}
Assumption \ref{Assumption: linearity} is similar to the assumption required in the path integral formulation of a single agent stochastic control problem \cite{satoh2016iterative}. A possible interpretation of condition \eqref{lambda} is that in a direction with high noise variance, the agent's control cost has to be low whereas that of the adversary has to be high. Therefore, the weights of the control costs $R_u$ and $R_v$ need to be tuned appropriately for the given diffusion coefficient $\Sigma(x, t)$ and the control gains $G_u(x, t)$ and $G_v(x, t)$ in the system dynamics \eqref{SDE}. See \cite{kappen2005path, williams2017model} for the further discussion on a similar condition in the single agent setting. Assumption \ref{Assumption: linearity} also implies that the stochastic noise has to enter the system dynamics via the control channels. 
Therefore, in what follows, we assume that system \eqref{SDE} can be partitioned into subsystems that are directly and non-directly driven by the noise as:
\begin{equation}\label{SDE partition}
    \begin{aligned}
  \begin{bmatrix}
   d\pmb{x}^{(1)} \\  d\pmb{x}^{(2)}
  \end{bmatrix}=& \begin{bmatrix}
    {f}^{(1)}(\pmb{x}, t) \\ {f}^{(2)}(\pmb{x}, t)
  \end{bmatrix}\!dt + \begin{bmatrix}
    \mathbf{0}\\{G_u}^{\!\!\!(2)}\!\left(\pmb{x}, t\right)
  \end{bmatrix}\!{u}(\pmb{x}, t)dt\\
  &\!\!\!\!\!\!+ \begin{bmatrix}
    \mathbf{0}\\{G_v}^{\!\!\!(2)}\!\left(\pmb{x}, t\right)
  \end{bmatrix}\!{v}(\pmb{x}, t)dt+\begin{bmatrix}
    \mathbf{0}\\{\Sigma}^{(2)}\!\left(\pmb{x}, t\right)
  \end{bmatrix}\!d\pmb{w}
\end{aligned}
\end{equation}
where $\mathbf{0}$ denotes a zero matrix of appropriate dimensions. By assuming a $\lambda$ satisfying Assumption \ref{Assumption: linearity} holds in (\ref{transformed HJB PDE}), we obtain the linear PDE in $\xi$ with Dirichlet boundary condition:
\begin{equation}\label{linearized risk-minimizing HJB}
 \!\!\begin{cases}
     \!\partial_t\xi\!=\!\frac{V\xi}{\lambda}\!-\!f^\top\partial_x\xi-\frac{1}{2}\text{Tr}\left(\Sigma\Sigma^\top\partial^2_x\xi\right),       & \forall(x,t)\!\in\!\mathcal{Q}, \\
    \!\!\underset{\substack{(x,t)\to(y,s) \\ (x,t)\in\mathcal{Q}}}{\lim}\xi(x,t)\!=\!\text{exp}{\left(-\frac{\phi(y)}{\lambda}\right)}, & \!\!\!\forall(y,s)\!\in\!\partial\mathcal{Q}.\\  
  \end{cases}
\end{equation}
The solution of a linear Dirichlet boundary value problem of the form (\ref{linearized risk-minimizing HJB}) exits under a sufficiently regular boundary condition, and it is unique \cite[Chapter 6]{friedman1975stochastic}. Furthermore, the solution admits the Feynman-Kac representation \cite{patil2022chance}. Suppose $\hat{\pmb{x}}(t)\in\mathbb{R}^n$ is an uncontrolled process driven by the following SDE:
\begin{equation}\label{uncontrolled SDE}
  d\hat{\pmb{x}}(t)\!=\!\!{f}\!\left(\hat{\pmb{x}}(t),\! t\right)\!dt\!+\!{\Sigma}\!\left(\hat{\pmb{x}}(t),\! t\right)\!d\pmb{w}(t)
\end{equation}
and let $ \hat{\pmb{t}}_{f} \coloneqq \text{inf}\{t> t_0: (\hat{\pmb{x}}(t), t)\notin \mathcal{Q}\}$. Then, the solution of the PDE (\ref{linearized risk-minimizing HJB}) is given as
\begin{equation}\label{xi}
\xi\left(x,t\right)=\mathbb{E}_{x, t}\left[\text{exp}{\left(-\frac{1}{\lambda}S\left(\tau\right)\right)}\right]  
\end{equation}
 where $S\left(\tau\right)$ denotes the cost-to-go of a trajectory $\tau$  of the uncontrolled system \eqref{uncontrolled SDE} starting at $(x,t)$:
\begin{equation}\label{Stau}
   S\left(\tau\right)=\phi\left(\hat{\pmb{x}}(\hat{\pmb{t}}_{f})\right)+\int_{t}^{\hat{\pmb{t}}_{f}} V\left(\hat{\pmb{x}}(t), t\right)dt. 
\end{equation}
Equation \eqref{xi} provides a path integral form for the exponentiated value function $\xi\left(x,t\right)$, which can be numerically evaluated using Monte Carlo sampling of trajectories generated by the uncontrolled SDE \eqref{uncontrolled SDE}.
We now obtain the expressions for the saddle-point policies via the following theorem:
\begin{theorem}
Suppose Assumption \ref{Assumption: linearity} holds and the system (\ref{SDE}) can be partitioned as \eqref{SDE partition}. Then, a saddle-point solution of the risk-minimizing zero-sum SDG \eqref{risk-minimizing SOC problem (tf)} exists, is unique and is given by
\begin{equation}\label{u* PI}
 \!\!\!\!u^*(x,t)dt\!=\!\mathcal{G}_u\!\left(x,t\right)\!\frac{\mathbb{E}_{x,t}\!\!\left[\text{exp}{\left(-\frac{1}{\lambda}S\left(\tau\right)\right)}\Sigma^{(2)}\!\!\left(x,t\right)d\pmb{w}\right]}{\mathbb{E}_{x,t}\left[\text{exp}{\left(-\frac{1}{\lambda}S\left(\tau\right)\right)}\right]}, 
\end{equation}
where 
\begin{equation*}
\mathcal{G}_u\!=\!R_u^{-1}{G_u^{(2)}}^\top\!\!\left(G_u^{(2)}R_u^{-1}{G_u^{(2)}}^\top\!\!\!-G_v^{(2)}R_v^{-1}{G_v^{(2)}}^\top\right)^{\!-1}
\end{equation*}
and 
\begin{equation}\label{v* PI}
 \!\!\!\!v^*(x,t)dt\!=\!\mathcal{G}_v\left(x,t\right)\!\frac{\mathbb{E}_{x,t}\!\!\left[\text{exp}{\left(-\frac{1}{\lambda}S\left(\tau\right)\right)}\Sigma^{(2)}\!\!\left(x,t\right)d\pmb{w}\right]}{\mathbb{E}_{x,t}\left[\text{exp}{\left(-\frac{1}{\lambda}S\left(\tau\right)\right)}\right]}, 
\end{equation}
where 
\begin{equation*}
\mathcal{G}_v\!=\!-R_v^{-1}{G_v^{(2)}}^\top\!\!\left(G_u^{(2)}R_u^{-1}{G_u^{(2)}}^\top\!\!\!-G_v^{(2)}R_v^{-1}{G_v^{(2)}}^\top\right)^{\!-1}.
\end{equation*}
\end{theorem}
\vspace{2mm}
\begin{proof}
The existence and uniqueness of the saddle-point solution follows from the existence and uniqueness of the linear Dirichlet boundary value problem \eqref{linearized risk-minimizing HJB} \cite[Chapter 6]{friedman1975stochastic} and from Theorem \ref{theorem: solution to risk-minimizing soc}. The saddle-point solution $u^*(x,t)$ (\ref{u*}) and $v^*(x,t)$ (\ref{v*}) can be computed by taking the gradient of (\ref{xi}) with respect to $x$ and using the condition (\ref{lambda}). (The derivation of (\ref{u*}) and (\ref{v*}) is in the same vein as the derivation of optimal controls in single agent settings \cite{satoh2016iterative,theodorou2010generalized}; not presented here for brevity.)
\end{proof}
 Equations (\ref{u* PI}) and (\ref{v* PI}) provide the path integral forms for the saddle-point equilibrium. Similar to \eqref{xi}, the expectations in (\ref{u* PI}) and (\ref{v* PI}) can be numerically evaluated in real-time via the Monte Carlo sampling of the trajectories generated by the uncontrolled SDE \eqref{uncontrolled SDE}. The path integral framework evaluates the solution locally without requiring knowledge of the solution nearby so that there is no need for a (global) discretization of the computational domain. This allows us to solve the game online without requiring any offline training or precomputations. Even though Monte Carlo simulations must be performed in real-time in order to evaluate (\ref{u* PI}, \ref{v* PI}) for the current $(x,t)$, these simulations can be massively parallelized through the use of GPUs. 

\section{Examples}\label{Sec: Examples}
In this section, we apply the path integral framework on two classes of risk-minimizing zero-sum SDGs (\ref{risk-minimizing SOC problem (tf)}): a disturbance attenuation problem and a pursuit-evasion game.

\subsection {Disturbance Attenuation Problem}\label{Sec: H-inf}
Consider a special class of systems (\ref{SDE}):
\begin{equation}\label{ex1:SDE}
\begin{aligned}
  d\pmb{x}\!=\!{f}(\pmb{x}, t)dt+{G_u}(\pmb{x}, t)\Big(\!{u}(\pmb{x}, t)dt \!+\! v\left(\pmb{x}, t\right)dt\!+\!d\pmb{w}\!\Big)
\end{aligned}
\end{equation}
where $u\left(\pmb{x}, t\right)\in\mathbb{R}^m$ is the control input, $v\left(\pmb{x},t\right)\in\mathbb{R}^m$ is the bounded disturbance and $\pmb{w}(t)\in\mathbb{R}^m$ is a Wiener process. Here, we have two sources of noise that corrupt the system's control input $u$: the bounded noise $v$ whose statistics are unknown and the white noise $d\pmb{w}$. 
In the disturbance attenuation problem, the objective is to design a policy $u$ in the presence of stochastic noise and bounded disturbance $v$ such that the system's control performance $ \mathbb{E}_{x_0, t_0}\!\left[\phi\left(\pmb{x}(\pmb{t}_{f})\right)\!+\!\int_{t_0}^{\pmb{t}_{f}}\!\left(\frac{1}{2}\pmb{u}^{\!\top}\!\pmb{u}+V\right)dt\right]$ is minimized. This problem can be solved using the following zero-sum SDG, where $u$ is considered as a control input of the first player (agent) and $v$ that of the second player (adversary):
\begin{equation}\label{J_gamma}
   \!\!\! \!\underset{u}{\min}\;\underset{v}{\max}\;\mathbb{E}_{x_0, t_0}\!\!\left[\!\phi\left(\pmb{x}(\pmb{t}_{f})\right)\!\!+\!\!\!\int_{t_0}^{\pmb{t}_{\!f}}\!\!\!\left(\!\frac{1}{2}\pmb{u}^\top\!\pmb{u}\!-\!\frac{\gamma^{2}}{\!2}\pmb{v}^{\!T}\!\pmb{v}\!+\!V\!\!\right)\!dt\!\right].
\end{equation}
$\gamma$ is a given positive constant which determines the level of disturbance attenuation. Theorem \ref{Propo. H_inf bound} provides an upper bound on the system's control performance (in the presence of a bounded disturbance $v$) that can be obtained by solving the game \eqref{J_gamma}.
\begin{theorem}\label{Propo. H_inf bound}
Suppose $(u^*_\gamma, v^*_\gamma)$ represent the saddle-point policies of the SDG (\ref{J_gamma}) for any $\gamma$, and let
\begin{equation}\label{delta gamma}
    \delta_\gamma \coloneqq \mathbb{E}_{x_0, t_0}^{u_\gamma^*, v_\gamma^*}\left[\int_{t_0}^{\pmb{t}_{\!f}}\!\!\!\!\pmb{v}^\top\!\pmb{v}\;dt\right]
\end{equation}
where the superscript on $\mathbb{E}$ denotes the polices under which the expectation is computed. Then, for all adversarial policies $v$ such that $\mathbb{E}_{x_0, t_0}^{u_\gamma^*, v}\left[\int_{t_0}^{\pmb{t}_{\!f}}\!\pmb{v}^\top\!\pmb{v}\;dt\right]\leq\delta$ (for any $\delta>0$), we get the following upper bound on the system's control performance in the presence of disturbance $v$:
\begin{equation}\label{UB on sys cost}
\begin{aligned}
      \mathbb{E}_{x_0, t_0}^{u_\gamma^*, v_\gamma^*}\!\!\left[\!\phi\left(\pmb{x}(\pmb{t}_{f})\right)\!\!+\!\!\!\int_{t_0}^{\pmb{t}_{\!f}}\!\!\!\left(\!\frac{1}{2}\pmb{u}^\top\!\pmb{u}\!+\!V\!\!\right)\!dt\!\right] + \frac{\gamma^2}{2}\left(\delta - \delta_\gamma\right)  \\
      \geq\mathbb{E}_{x_0, t_0}^{u_\gamma^*, v}\!\!\left[\!\phi\left(\pmb{x}(\pmb{t}_{f})\right)\!+\!\!\!\int_{t_0}^{\pmb{t}_{\!f}}\!\!\!\left(\!\frac{1}{2}\pmb{u}^\top\!\pmb{u}+\!V\!\!\right)\!dt\!\right].
\end{aligned}
\end{equation}
\end{theorem}
\vspace{1mm}
\begin{proof}
See Appendix B.
\end{proof}
\vspace{2mm}
In order to solve the HJI PDE associated with the game (\ref{J_gamma}) via the path integral framework described in Section \ref{Sec: PI}, it is necessary to find a constant $\lambda>0$ (by Assumption \ref{Assumption: linearity}) such that
\begin{equation*}\label{linearity on H_inf}
    \lambda\left(1-\frac{1}{\gamma^2}\right)=1.
\end{equation*}
Therefore, for all $\gamma>1$, Assumption 1 is satisfied and as a consequence, the zero-sum SDG (\ref{J_gamma}) admits a unique saddle-point solution. \par
We now present a simulation study of the disturbance attenuation problem using a unicycle navigation example. Consider the following unicycle dynamics model:
\begin{equation} \label{unicycle model}
\begin{aligned}
    \begin{bmatrix}
    d\pmb{p}_x\\d\pmb{p}_y\\d\pmb{s}\\d\pmb{\theta}
    \end{bmatrix}\!=&\!
    -k
    \begin{bmatrix}
    \pmb{p}_x\\
    \pmb{p}_y\\
    \pmb{s}\\
    \pmb{\theta}
    \end{bmatrix}dt+
    \begin{bmatrix}
    \pmb{s}\cos{\pmb{\theta}}\\\pmb{s}\sin{\pmb{\theta}}\\0\\0
    \end{bmatrix}\!dt\!\\
    &\!\!\!\!+\begin{bmatrix}
    0 & 0\\0 & 0\\1 & 0\\0 & 1
    \end{bmatrix} \!
    \left(\begin{bmatrix}
    a\\
    \omega
    \end{bmatrix}\!dt \!+\!
    \begin{bmatrix}
    \Delta a\\
    \Delta \omega
    \end{bmatrix}\!dt \!+\!
    \begin{bmatrix}
    \sigma & 0\\
    0 & \nu 
    \end{bmatrix}\!d\pmb{w}
    \right),
\end{aligned}
\end{equation}
where $(\pmb{p}_x,\:\pmb{p}_y)$, $\pmb{s}$ and $\pmb{\theta}$ denote the position, speed, and the heading angle of the unicycle, respectively. The control input $u\coloneqq\begin{bmatrix} a & \omega \end{bmatrix}^\top$ consists of acceleration $a$ and angular speed $\omega$. $v\coloneqq\begin{bmatrix} \Delta a & \Delta \omega \end{bmatrix}^\top$ is the bounded disturbance acting on the system's control input, and $d\pmb{w}\in\mathbb{R}^2$ is the white noise with $\sigma$ and $\nu$ being the noise level parameters. As illustrated in Figure \ref{Fig. effect of eta}, the unicycle aims to navigate in a two-dimensional space from its initial position (represented by the yellow star) to the origin (represented by the magenta star), in finite time, while avoiding the red obstacles and the outer boundary. The white region that lies between the outer boundary and the obstacles is the safe region $\mathcal{X}_s$. This is a disturbance attenuation problem, since the unicycle aims to design its control policy $u$ in order to minimize the control performance and risk of failure (collision with the obstacles or the outer boundary) under worst-case disturbance $v$. Therefore, we can formulate this problem as the risk-minimizing zero-sum SDG (\ref{J_gamma}). In the simulation, we set $\sigma=\nu=0.1$, $k=0.2$, $t_0=0$, $T=10$, $x_0=\begin{bmatrix}
-0.4 &-0.4 & 0 &0
\end{bmatrix}^\top$, $V(\pmb{x}) = \pmb{p}_x^2 + \pmb{p}_y^2$ and $\psi\left(\pmb{x}(T)\right) = \pmb{p}_x^2(T) + \pmb{p}_y^2(T)$. In order to evaluate the optimal policies (\ref{u* PI}) and (\ref{v* PI}) via Monte Carlo sampling, $10^4$ trajectories and a step size equal to $0.01$ are used. We demonstrate two experiments. 
\subsubsection{Experiment 1}
\begin{figure}
     \centering
       \begin{tabular}{c c}
       \!\!\!\!\!\!\!\!\!\includegraphics[scale=0.35]{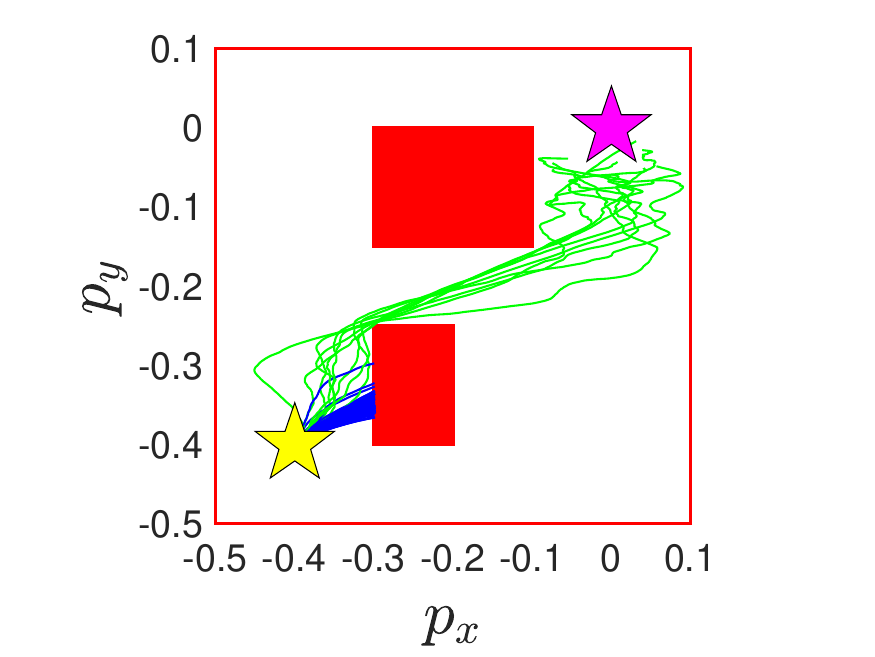} &\!\!\!\!\!\!\!\!\!\!\!\!\!\!\!\!\!\!\!\!\!\!\!\!\includegraphics[scale=0.35]{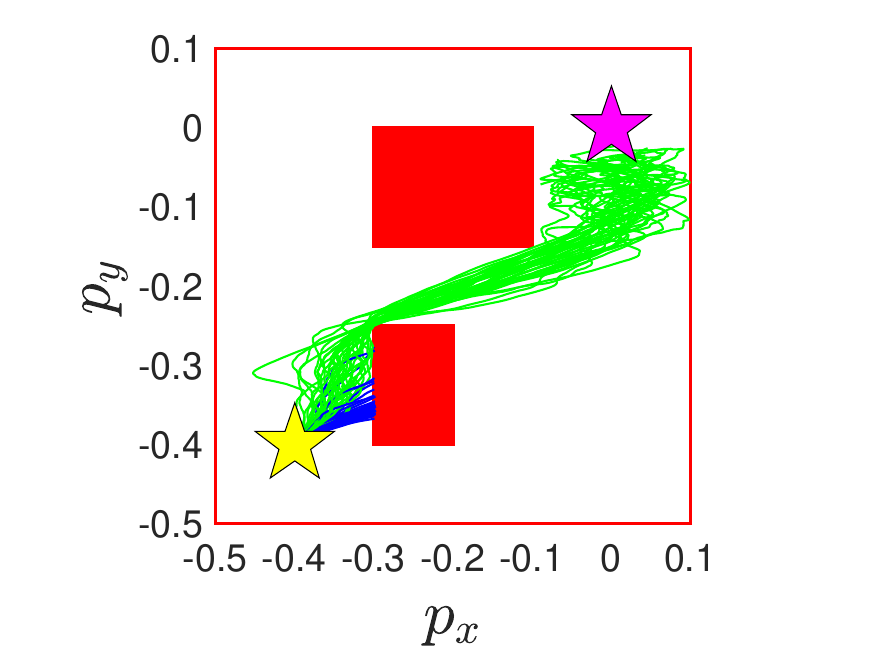} \\
       \!\!\!\!\!\!\!\!\!(a) $\gamma^2 = 2$, $P^\mathrm{ag}_{\mathrm{fail}} = 0.9$  & \!\!\!\!\!\!\!\!\!\!\!\!\!\!\!\!\!\!\!\!\!\!\!\!\!(b) $\gamma^2 = 7 $, $P^\mathrm{ag}_{\mathrm{fail}}= 0.64$ \\

       \end{tabular}
         \caption{Unicycle navigation in the presence of bounded and stochastic disturbances. The start position is shown by a yellow star and the target position (the origin) by a magenta star. $100$ sample trajectories generated using saddle-point policies ($u^*, v^*$) for two values of $\gamma$ are shown. The trajectories are color-coded; blue paths collide with the red obstacles or the outer boundary, while the green paths converge in the neighborhood of the magenta star. The failure probabilities of the agent $P^\mathrm{ag}_{\mathrm{fail}}$ are noted below each case.} 
         \label{Fig. effect of eta}
 \end{figure}

In this experiment, we set $\eta = 0.67$ and plot in Figure \ref{Fig. effect of eta} $100$ sample trajectories generated using synthesized saddle-point policies $(u^*, v^*)$ for two values of $\gamma$. The trajectories are color-coded; the blue paths collide with the obstacles, while the green paths converge in the neighborhood of the origin (the target position). The figure shows that for a higher value of $\gamma$ i.e., when the adversary becomes less powerful, the failure probability of the agent $P^\mathrm{ag}_{\mathrm{fail}}$ reduces. 
\subsubsection{Experiment 2}
\begin{figure}
     \centering
       \begin{tabular}{c c}
       \!\!\!\!\!\!\!\!\!\includegraphics[scale=0.35]{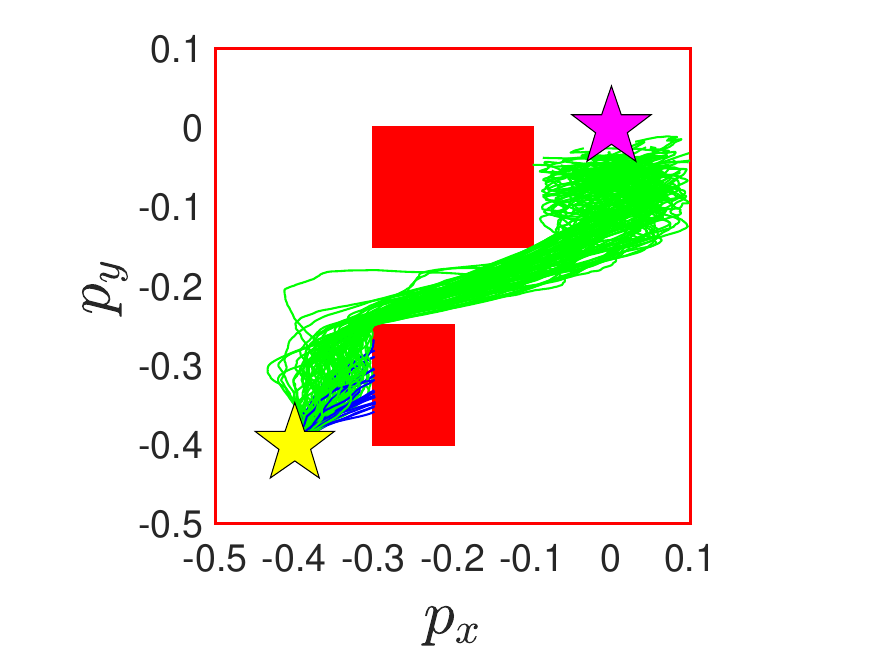} &\!\!\!\!\!\!\!\!\!\!\!\!\!\!\!\!\!\!\!\!\!\!\!\!\!\includegraphics[scale=0.35]{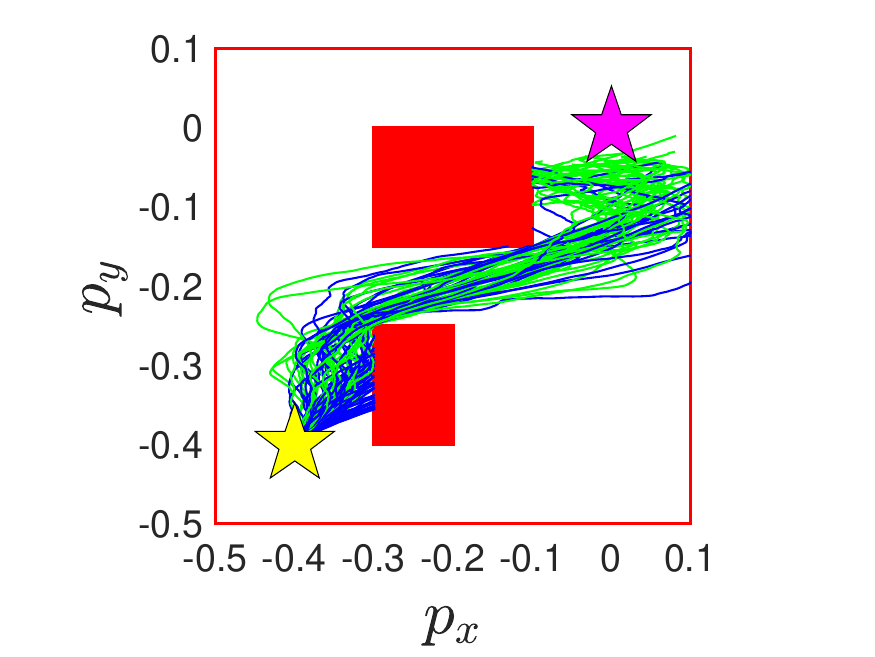} \\
       \!\!\!\!\!\!\!\!\!(a) Agent is aware of the& \!\!\!\!\!\!\!\!\!\!\!\!\!\!\!\!\!\!(b) Agent is not aware of \\
       \!\!\!\!\!\!\!\!\!adversary, $P^\mathrm{ag}_{\mathrm{fail}} = 0.23 $ & \!\!\!\!\!\!\!\!\!\!\!\!\!\!\!\!\!\!the adversary,  $P^\mathrm{ag}_{\mathrm{fail}}= 0.65$
      
       \end{tabular}
         \caption{Unicycle navigation in the presence of bounded and stochastic disturbances. (a) Agent is aware of the presence of adversary. $100$ sample trajectories generated using saddle-point policies $(u^*, v^*)$. (b) Agent is not aware of the presence of adversary. $100$ sample trajectories generated using $({\widetilde{u}}^*, v^*)$. The failure probabilities of the agent $P^\mathrm{ag}_{\mathrm{fail}}$ are noted for each case.} 
         \label{Fig. safe unsafe}
 \end{figure}

In this experiment, we set $\gamma^2=3$, $\eta=1$ and study the effect of ignoring the adversary. First, we compute saddle-point policies $(u^*, v^*)$ for the game (\ref{J_gamma}) same as Experiment 1 and plot in Figure \ref{Fig. safe unsafe}-(a) $100$ sample trajectories generated using $(u^*, v^*)$. In this case, the agent is aware of the adversary and designs its policy $u^*$ cautiously. The probability of failure is $23\%$. In the second case, the agent is not aware of the presence of adversary and computes its policy (say) ${\widetilde{u}}^*$ by solving a single agent optimization problem. However, in reality the adversary is present and suppose it follows the policy $v^*$. Figure \ref{Fig. safe unsafe}-(b) shows $100$ sample trajectories generated using $({\widetilde{u}}^*, v^*)$. In this case, the agent's performance is poor, it fails $65\%$ of the times. The color-coding of the trajectories is same as Experiment 1.
      

\subsection {Pursuit-Evasion Game}\label{Sec: PE}
Consider a two-player zero-sum SDG on a finite time horizon $[t_0, T]$, in which the adversary is chasing the agent and the agent is trying to escape from the adversary. We will call the adversary as a pursuer and the agent as an evader. Suppose the evader and the pursuer are moving in a two-dimensional plane according to
 \begin{equation}\label{sim EP model}
     \begin{split}
    d\pmb{p}^E_x={u}_xdt+\sigma^E_x d\pmb{w}^E_x,\qquad d\pmb{p}^P_x={v}_xdt+\sigma^P_x d\pmb{w}^P_x,\\
     d\pmb{p}^E_y={u}_ydt+\sigma^E_y d\pmb{w}^E_y,\qquad d\pmb{p}^P_y={v}_ydt+\sigma^P_y d\pmb{w}^P_y,
\end{split}
 \end{equation}
where $\pmb{x}_E\coloneqq\begin{bmatrix}
 \pmb{p}^E_x & \pmb{p}^E_y 
\end{bmatrix}^\top$ is the position and $u \coloneqq\begin{bmatrix}
 u_x & u_y 
\end{bmatrix}^\top$ is the control input of the evader. Similarly,  $\pmb{x}_P\coloneqq\begin{bmatrix}\pmb{p}^P_x & \pmb{p}^P_y\end{bmatrix}^\top$ and $v\coloneqq\begin{bmatrix}
 v_x &v_y 
\end{bmatrix}^\top$ are the position and control input of the pursuer. $\pmb{w}_x^E, \pmb{w}_y^E, \pmb{w}_x^P, \pmb{w}_y^P$ are independent one-dimensional standard Brownian motions. If at any time $t \in (t_0, T]$, the pursuer gets within a distance $\rho$ of the evader, then it catches the evader and the evader fails. On the other hand, if the evader avoids getting within a distance $\rho$ of the pursuer for the entire time horizon $[t_0, T]$, then that's a failure for the pursuer. The pursuer aims at designing its control policy $v$ in order to maximize the probability of catching the evader, whereas the evader seeks the opposite by designing $u$. For this two-player differential game, it is the relative position of the pursuer and evader that is important (and relevant), rather than their absolute positions. Let $\pmb{x}\coloneqq\begin{bmatrix}
  \pmb{p}_x & \pmb{p}_y
\end{bmatrix}$ be the evader's position with respect to the pursuer where 
\begin{equation*}
    \pmb{p}_x = \pmb{p}^E_x - \pmb{p}^P_x, \quad \pmb{p}_y = \pmb{p}^E_y - \pmb{p}^P_y
\end{equation*}
and the origin coincides with the pursuer's position. Thus, the coordinate system is attached to the pursuer and is not fixed in space. The system $\pmb{x}$ follows the SDE
\begin{equation}\label{rel. dynamics}
    d\pmb{x}=\begin{bmatrix}
    d\pmb{p}_x\\ d\pmb{p}_y
    \end{bmatrix} = 
    \begin{bmatrix}
    u_x\\ u_y
    \end{bmatrix}dt - 
    \begin{bmatrix}
    v_x\\ v_y
    \end{bmatrix}dt +
    \begin{bmatrix}
  \sigma_x & 0\\  0 & \sigma_y
    \end{bmatrix}d\pmb{w}
\end{equation}
where $\sigma_x =\sqrt{(\sigma_x^E)^2 + (\sigma_x^P)^2}$, $\sigma_y =\sqrt{(\sigma_y^E)^2 + (\sigma_y^P)^2}$ and $\pmb{w}$ is a two-dimensional standard Brownian motion. In this game, the safe set $\mathcal{X}_s$ can be defined as $\mathcal{X}_s \coloneqq \left\{x\in\mathbb{R}^2:\|x\|>\rho\right\}$. Suppose the control cost matrix $R_u$ of the evader is unity and that of the pursuer $R_v = {r_v}^2$, where $r_v$ is a given positive scalar constant. Therefore, the risk-minimizing zero-sum SDG takes the form:
\begin{equation}\label{SDG of PE}
    \!\!\underset{u}{\min}\;\underset{v}{\max}\;\mathbb{E}_{x_0, t_0}\!\!\left[\!\phi\left(\pmb{x}(\pmb{t}_{f}\!)\right)\!+\!\!\!\int_{t_0}^{\pmb{t}_{\!f}}\!\!\!\!\left(\!\frac{1}{2}\pmb{u}^{\!T}\!\pmb{u}\!-\!\frac{{r_{\!v}}^{\!2}}{2}\pmb{v}^{\!T}\!\pmb{v}\!+\!V\!\!\right)\!dt\!\right]\!.
\end{equation}
In order to solve the associated HJI equation of this game via the path integral framework, it is necessary to find a constant $\lambda>0$ (by Assumption 1) such that
\begin{equation*}
    \lambda\left(1-\frac{1}{{r_v}^2}\right)=1.
\end{equation*}
Therefore, for all $r_v>1$, Assumption 1 is satisfied and as a consequence, the zero-sum SDG (\ref{SDG of PE}) admits a unique saddle-point solution. \par 

\begin{figure}
     \centering
       \begin{tabular}{c}
       \!\!\!\!\!\!\!\!\!\!\!\!\!\!\!\!\!\includegraphics[scale=0.17]{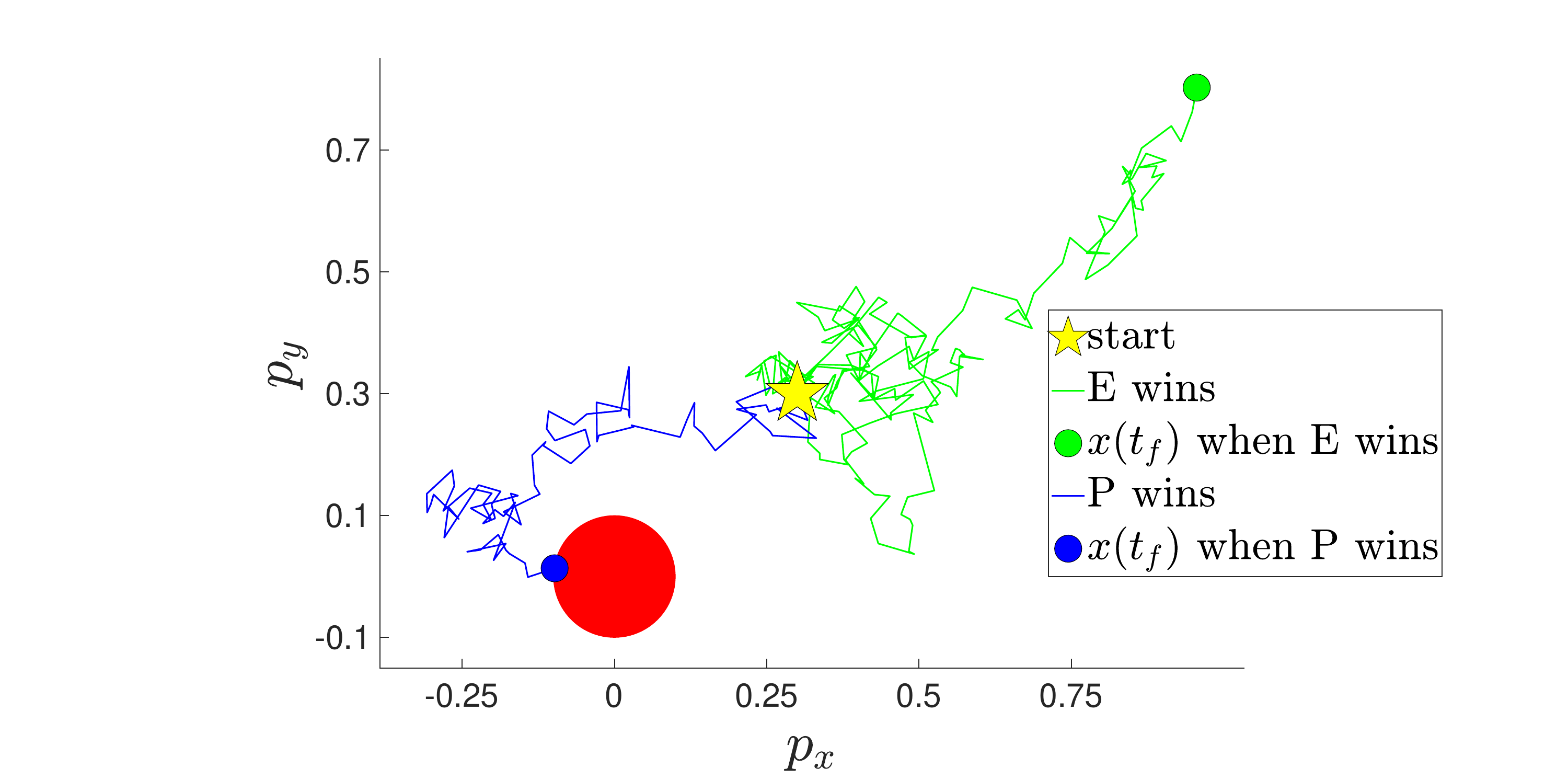} \\
       \end{tabular}
         \caption{Two sample trajectories of the relative position of the players in a pursuit-evasion game. The start position of the trajectories is shown by a yellow star. The red disc of radius $\rho = 0.1$, centered at the origin represents that the pursuer is within the distance $\rho$ of the evader. The green trajectory never enters the red disc in the horizon $[t_0, T]$, thus, it represents a case when the evader wins. The blue trajectory, enters the red disc and thus represents a case when the pursuer wins.} 
         \label{Fig. EP}
 \end{figure}

\begin{figure}
    \centering
\includegraphics[scale=0.15]{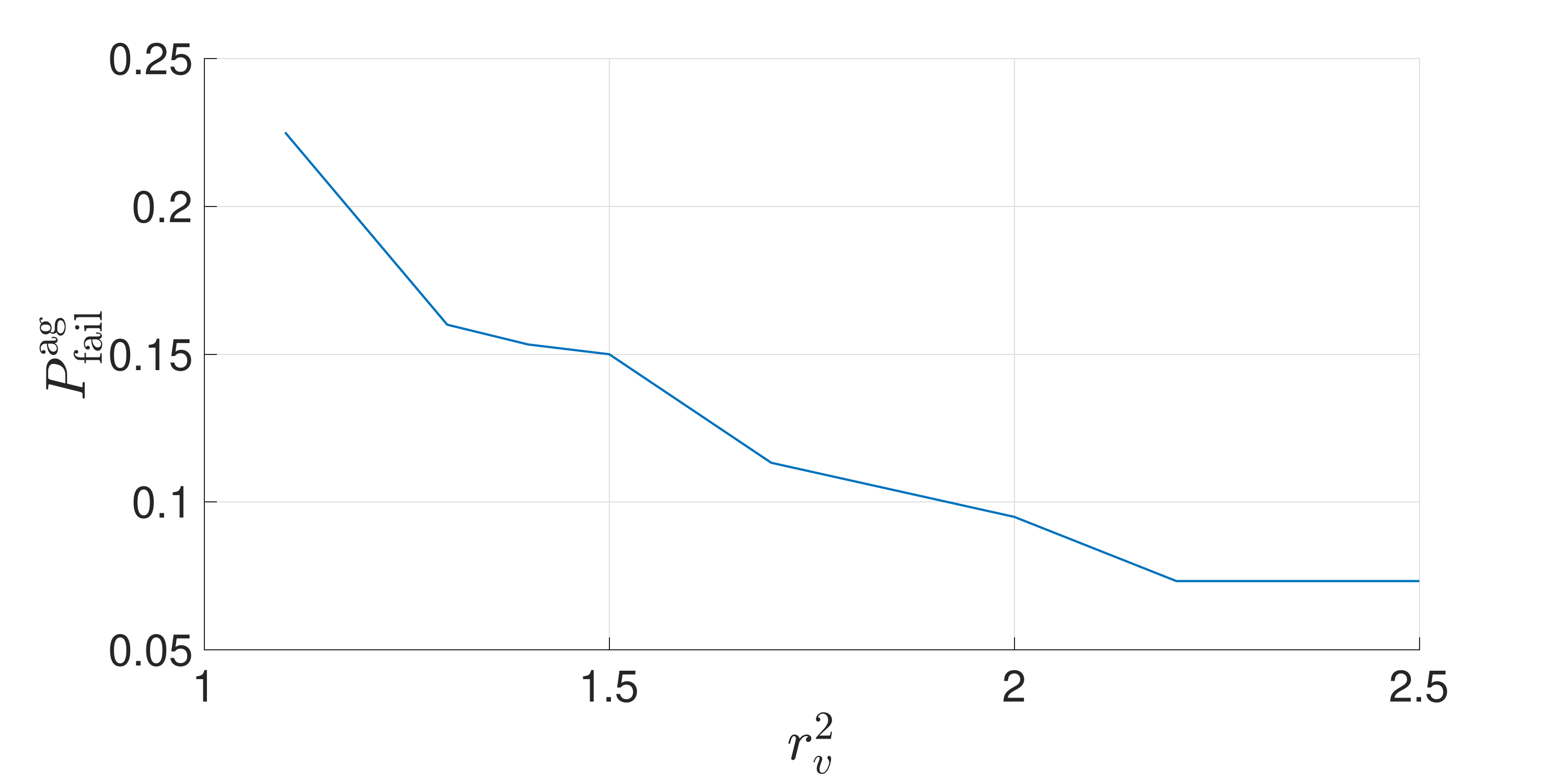} 
        \caption{Failure probabilities of the agent (i.e., evader) as a function of $r_v$, when the players follow the saddle-point policies ($u^*, v^*$).} 
        \label{Fig. pfail with rv2}
\end{figure}

In the simulation, we set $\sigma_x^E=\sigma_y^E=\sigma_x^P=\sigma_y^P=\sqrt{0.1}$, $\rho=0.1$, $t_0=0$, $T=2$, $x_0 = \begin{bmatrix}
  0.3 & 0.3
\end{bmatrix}^\top$, $V(\pmb{x})=\psi\left(\pmb{x}(T)\right)=0$,  $\eta=0.2$, ${r_v}^2 = 2$. Figure \ref{Fig. EP} shows a plot of two sample trajectories of system (\ref{rel. dynamics}) generated using synthesized saddle-point policies $(u^*, v^*)$. The trajectories start from $x_0$ shown by the yellow star. The red disc of radius $\rho = 0.1$, centered at the origin represents that the pursuer is within a distance $\rho$ of the evader. The green trajectory never enters the red disc in the horizon $[t_0, T]$, thus, it represents a case when the evader escapes from the pursuer. The blue trajectory on the other hand, enters the red disc and thus represents a case when the pursuer catches the evader. Figure \ref{Fig. pfail with rv2} shows a plot of failure probabilities of the agent (i.e., evader) as a function of $r_v$, when the players follow the saddle-point policies ($u^*, v^*$). These values are computed using na\"ive Monte Carlo sampling, with $400$ sample trajectories. The plot shows that as the control cost weight $r_v$ of the adversary (i.e., pursuer) increases, the chances of evader getting caught reduces.   

 \par
The aim of the presented simulation studies is to validate the proposed theoretical formulation of the risk-minimizing zero-sum SDGs. Future work will emphasize on scaling this framework to higher dimensional and more complex game dynamics.

\section{Conclusion} \label{sec: Conclusion}
The paper presented an HJI-PDE-based solution approach for a risk-minimizing two-player zero-sum stochastic differential game (SDG). This is a variable-end-time game in which each player tries to balance the trade-off between the probability of failure and the control cost. A sufficient condition for a saddle-point solution of this game was derived and it was shown that this game can be solved via an HJI PDE with the Dirichlet boundary condition. We developed a path integral framework to numerically solve a class of risk-minimizing zero-sum SDGs whose associated HJI PDE can be linearized and established the existence and uniqueness of the saddle-point solution. The presented approach allows the game to be solved online without the need of any offline training or precomputations. Application of our approach on two classes of risk-minimizing zero-sum SDGs: a disturbance attenuation problem and a pursuit-evasion game was presented and the framework was validated through simulation studies.\par 

In the future, we plan to conduct sample complexity analysis for path integral control in order to investigate how the accuracy of Monte Carlo sampling affects the solution of SDGs. The central challenge in using the path integral framework is the particular requirement on the relationship between the cost function and the noise covariance. This requirement restricts the class of applicable system models and cost functions. In the future work, we plan to find alternatives in order to get rid of this restrictive requirement (one such solution is provided in \cite{satoh2016iterative}). Another topic of future investigation could be chance-constrained stochastic games in which each player would aim to satisfy a hard bound on its failure probability.

\appendix
\subsection{Proof of Theorem 1}
Let $J(x,t)$ be the function satisfying (a) and (b). By Dynkin's formula \cite{oksendal2013stochastic,durrett2019probability}, for each $(x,t)\in\overline{\mathcal{Q}}$ we have
\begin{align}\label{ito}
   \mathbb{E}_{x,t}\left[J\left(\pmb{x}(\pmb{t}_{f}), \pmb{t}_{f}\right)\right]&=J(x,t)\\ &\!\!\!\!\!\!\!\!\!\!\!\!\!\!\!\!\!\!\!\!\!\!\!\!\!\!\!\!\!\!\!\!\!\!\!\!\!\!\!\!\!\!\!\!\!\!\!\!\!\!\!\!+\!\mathbb{E}_{x,t}\!\!\left[\!\int_{t}^{\pmb{t}_{\!f}}\!\!\!\!\left(\!\!\partial_tJ \!+\! (f\!+\!G_{\!u}u\!+\!G_{\!v}v)^{\!T}\!(\partial_xJ)\!+\!\frac{1}{2}\text{Tr}\!\left(\Sigma\Sigma^\top\!\partial_x^2J\right)\!\!\right)\!ds\!\right]\!\!.\notag
\end{align}
By the boundary condition of the PDE (\ref{HJB PDE}),\newline $J\left(\pmb{x}(\pmb{t}_{f}), \pmb{t}_{f}\right) =\phi\left(\pmb{x}(\pmb{t}_{f})\right)$. Hence, from (\ref{ito}), we obtain
\begin{align}\label{ito after plugging BC}
   J(x,t)&=\mathbb{E}_{x,t}\left[\phi\left(\pmb{x}(\pmb{t}_{f})\right)\right]\\ &\!\!\!\!\!\!\!\!\!\!\!\!\!\!\!\!\!\!\!\!\!\!\!-\!\mathbb{E}_{x,t}\!\!\left[\!\int_{t}^{\pmb{t}_{\!f}}\!\!\!\!\left(\!\!\partial_tJ \!+\! (f\!+\!G_{\!u}u\!+\!G_{\!v}v)^{\!T}\!(\partial_{x}J)\!+\!\frac{1}{2}\text{Tr}\!\left(\Sigma\Sigma^{T}\!\partial_x^2J\right)\!\!\right)\!ds\!\right]\!\!.\notag
\end{align}
Now, notice that the right hand side of the PDE in (\ref{HJB PDE}) can be expressed as the minimum and maximum value of a quadratic form in $u$ and $v$, respectively, as follows:
\begin{equation}\label{-dJ/dt as minimax}
    \begin{aligned}
         \!\!-\partial_tJ\!=\!\min_{{u}} \max_{{v}}\!\Bigg[&\frac{1}{2}u^\top\!R_uu-\frac{1}{2}v^\top\!R_vv\!+\!V\!\\
         &\!\!\!\!\!\!\!\!\!\!\!\!\!\!\!\!\!\!\!\!\!\!\!\!+\!\left(\!f\!+\!G_uu\!+\!G_vv\right)^\top\!\partial_xJ+\!\frac{1}{2}\text{Tr}\!\left(\Sigma\Sigma^\top\partial^2_xJ\right)\!\Bigg]\!.
          \end{aligned} 
\end{equation}
Observe that the ``$\underset{u}{\min}$", ``$\underset{v}{\max}$" operations in (\ref{-dJ/dt as minimax}) can be interchanged. Hence, the game formulated in (\ref{risk-minimizing SOC problem (tf)}) satisfies the \textit{Isaacs condition} \cite{isaacs1999differential}. If $\hat{u}$ and $\hat{v}$ represent the minimum and maximum values of the right hand side of (\ref{-dJ/dt as minimax}) respectively, then 
\begin{equation}\label{uhat, vhat}
    \hat{u} = -R_u^{-1}G_u^\top\partial_xJ, \qquad \hat{v} =R_v^{-1}G_v^\top\partial_xJ.
\end{equation}
Therefore, for an arbitrary $u$, we have 
\begin{equation}\label{pde ineqaulity}
    \begin{aligned}
         \!\!\!\!-\partial_tJ\!\leq\!\Bigg[&\frac{1}{2}u^\top\!R_uu-\frac{1}{2}{\hat{v}}^\top R_v\hat{v}\!+\!V\!\\
         &+\!\left(\!f\!+\!G_uu\!+\!G_v\hat{v}\right)^\top\!\partial_xJ+\!\frac{1}{2}\text{Tr}\!\left(\Sigma\Sigma^\top\partial^2_xJ\right)\!\Bigg]\!.
          \end{aligned} 
\end{equation}
Now, notice that the equality in (\ref{ito after plugging BC}) holds for any $v$. Replacing $v$ by $\hat{v}$ in (\ref{ito after plugging BC}) yields
\begin{equation}\label{ito after plugging BC2}
    \begin{aligned}
       J(x,t)&=\mathbb{E}_{x,t}\left[\phi\left(\pmb{x}(\pmb{t}_{f})\right)\right]\\ &\!\!\!\!\!\!\!\!\!\!\!\!\!\!\!\!\!\!\!\!\!\!\!-\!\mathbb{E}_{x,t}\!\!\left[\!\int_{t}^{\pmb{t}_{\!f}}\!\!\!\!\left(\!\!\partial_tJ \!+\! (\!f\!+\!G_{\!u}u\!+\!G_{\!v}\hat{v})^{\!T}\!(\partial_xJ)\!+\!\frac{1}{2}\text{Tr}\!\left(\Sigma\Sigma^\top\!\partial_x^2J\right)\!\!\right)\!\!ds\!\right]\!\!.
    \end{aligned}
\end{equation}
Combining (\ref{pde ineqaulity}) and (\ref{ito after plugging BC2}), we obtain

\begin{equation}\label{J leq C_hat}
  \begin{aligned}
       \!\!\!\!\!\!J(x,t)\!&\leq\!\mathbb{E}_{x,t}\!\!\left[\!\phi\!\left(\pmb{x}(\pmb{t}_{f}\!)\right)\!+\!\!\!\int_{t}^{\pmb{t}_{f}}\!\!\!\!\left(\!\frac{1}{2}\pmb{u}^{T\!}\!R_u\pmb{u}\!-\!\frac{1}{2}\hat{\pmb{v}}^{\!T}\!\!R_v\hat{\pmb{v}}+\!V\!\!\right)\!ds\!\right]\\
       &= C\left(x,t;u,\hat{v}\right)
    \end{aligned}   
\end{equation}
where the equality holds iff $ \hat{u} =-R_u^{-1}G_u^\top\partial_xJ$. Similarly, for an arbitrary $v$, we can show that 
\begin{equation}\label{J geq C_hat}
  J(x,t)\geq C\left(x,t;\hat{u},v\right)
\end{equation}
where the equality holds iff $\hat{v} =R_v^{-1}G_v^\top\partial_xJ$.
Therefore, from Definition \ref{Def: SP}, it follows that the pair of policies $(\hat{u}, \hat{v})$ defined in (\ref{uhat, vhat}) provides the optimal solution to the zero-sum game formulated in Problem \ref{Problem: Risk-minimizing SOC problem} and $J(x,t)$ is the value of the game.
\subsection{Proof of Theorem 3}
Consider cost of the SDG (\ref{J_gamma}) under the saddle-point policies $(u^*_\gamma, v^*_\gamma)$:

\begin{subequations}\label{contra. final}
\begin{align}
&\mathbb{E}_{x_0, t_0}^{u_\gamma^*, v_\gamma^*}\!\!\left[\!\phi\left(\pmb{x}(\pmb{t}_{f})\right)\!\!+\!\!\!\int_{t_0}^{\pmb{t}_{\!f}}\!\!\!\left(\!\frac{1}{2}\pmb{u}^\top\pmb{u}\!+\!V\!\!-\!\frac{\gamma^2}{2}\pmb{v}^\top\pmb{v}\!\!\right)\!dt\!\right] \nonumber\\
=&\;\mathbb{E}_{x_0, t_0}^{u_\gamma^*, v_\gamma^*}\!\!\left[\!\phi\left(\pmb{x}(\pmb{t}_{f})\right)\!\!+\!\!\!\int_{t_0}^{\pmb{t}_{\!f}}\!\!\!\left(\!\frac{1}{2}\pmb{u}^\top\pmb{u}\!+\!V\!\!\right)\!dt\!\right] -\frac{\gamma^2}{2}\delta_\gamma \label{contra. final1}\\
\geq&\;\mathbb{E}_{x_0, t_0}^{u_\gamma^*, v}\!\!\left[\!\phi\left(\pmb{x}(\pmb{t}_{f})\right)\!\!+\!\!\!\int_{t_0}^{\pmb{t}_{\!f}}\!\!\!\left(\!\frac{1}{2}\pmb{u}^\top\pmb{u}\!+\!V\!\!-\!\frac{\gamma^2}{2}\pmb{v}^\top\pmb{v}\!\!\right)\!dt\!\right]\label{contra. final2}\\
\geq&\;\mathbb{E}_{x_0, t_0}^{u_\gamma^*, v}\!\!\left[\!\phi\left(\pmb{x}(\pmb{t}_{f})\right)\!\!+\!\!\!\int_{t_0}^{\pmb{t}_{\!f}}\!\!\!\left(\!\frac{1}{2}\pmb{u}^\top\pmb{u}\!+\!V\!\!\right)\!dt\!\right] -\frac{\gamma^2}{2}\delta.\label{contra. final3}
\end{align}
\end{subequations}
The equation \eqref{contra. final1} follows from (\ref{delta gamma}). For any adversarial policy $v$, the inequality \eqref{contra. final2} follows because $v^*_\gamma$ maximizes the cost in \eqref{J_gamma}. The inequality \eqref{contra. final3} follows from the bound on $\mathbb{E}_{x_0, t_0}^{u_\gamma^*, v}\left[\int_{t_0}^{\pmb{t}_{\!f}}\!\pmb{v}^\top\!\pmb{v}\;dt\right]$. Using \eqref{contra. final1} and \eqref{contra. final3}, we get the desired inequality \eqref{UB on sys cost}. 
\bibliographystyle{IEEEtran}
\bibliography{root}

\begin{thebibliography}{10}
\providecommand{\url}[1]{#1}
\csname url@rmstyle\endcsname
\providecommand{\newblock}{\relax}
\providecommand{\bibinfo}[2]{#2}
\providecommand\BIBentrySTDinterwordspacing{\spaceskip=0pt\relax}
\providecommand\BIBentryALTinterwordstretchfactor{4}
\providecommand\BIBentryALTinterwordspacing{\spaceskip=\fontdimen2\font plus
\BIBentryALTinterwordstretchfactor\fontdimen3\font minus
  \fontdimen4\font\relax}
\providecommand\BIBforeignlanguage[2]{{%
\expandafter\ifx\csname l@#1\endcsname\relax
\typeout{** WARNING: IEEEtran.bst: No hyphenation pattern has been}%
\typeout{** loaded for the language `#1'. Using the pattern for}%
\typeout{** the default language instead.}%
\else
\language=\csname l@#1\endcsname
\fi
#2}}

\bibitem{bacsar1998dynamic}
T.~Ba{\c{s}}ar and G.~J. Olsder, \emph{Dynamic noncooperative game
  theory}.\hskip 1em plus 0.5em minus 0.4em\relax SIAM, 1998.

\bibitem{nahin2012chases}
P.~J. Nahin, ``Chases and escapes,'' in \emph{Chases and Escapes}.\hskip 1em
  plus 0.5em minus 0.4em\relax Princeton University Press, 2012.

\bibitem{7172219}
W.~Sun and P.~Tsiotras, ``Pursuit evasion game of two players under an external
  flow field,'' in \emph{2015 American Control Conference (ACC)}, 2015, pp.
  5617--5622.

\bibitem{falcone2006numerical}
M.~Falcone, ``Numerical methods for differential games based on partial
  differential equations,'' \emph{International Game Theory Review}, vol.~8,
  no.~02, pp. 231--272, 2006.

\bibitem{huang2014automation}
H.~Huang, J.~Ding, W.~Zhang, and C.~J. Tomlin, ``Automation-assisted
  capture-the-flag: A differential game approach,'' \emph{IEEE Transactions on
  Control Systems Technology}, vol.~23, no.~3, pp. 1014--1028, 2014.

\bibitem{mitchell2005time}
I.~M. Mitchell, A.~M. Bayen, and C.~J. Tomlin, ``A time-dependent
  hamilton-jacobi formulation of reachable sets for continuous dynamic games,''
  \emph{IEEE Transactions on automatic control}, vol.~50, no.~7, pp. 947--957,
  2005.

\bibitem{vrabie2011adaptive}
D.~Vrabie and F.~Lewis, ``Adaptive dynamic programming for online solution of a
  zero-sum differential game,'' \emph{Journal of Control Theory and
  Applications}, vol.~9, no.~3, pp. 353--360, 2011.

\bibitem{prajapat2021competitive}
M.~Prajapat, K.~Azizzadenesheli, A.~Liniger, Y.~Yue, and A.~Anandkumar,
  ``Competitive policy optimization,'' in \emph{Uncertainty in Artificial
  Intelligence}.\hskip 1em plus 0.5em minus 0.4em\relax PMLR, 2021, pp. 64--74.

\bibitem{liu2020adaptive}
M.~Liu, Y.~Wan, F.~L. Lewis, and V.~G. Lopez, ``Adaptive optimal control for
  stochastic multiplayer differential games using on-policy and off-policy
  reinforcement learning,'' \emph{IEEE transactions on neural networks and
  learning systems}, vol.~31, no.~12, pp. 5522--5533, 2020.

\bibitem{lin2017multiagent}
X.~Lin, P.~A. Beling, and R.~Cogill, ``Multiagent inverse reinforcement
  learning for two-person zero-sum games,'' \emph{IEEE Transactions on Games},
  vol.~10, no.~1, pp. 56--68, 2017.

\bibitem{artzner1999coherent}
P.~Artzner, F.~Delbaen, J.-M. Eber, and D.~Heath, ``Coherent measures of
  risk,'' \emph{Mathematical finance}, vol.~9, no.~3, pp. 203--228, 1999.

\bibitem{dixit2022risk}
A.~Dixit, M.~Ahmadi, and J.~W. Burdick, ``Risk-averse receding horizon motion
  planning,'' \emph{arXiv preprint arXiv:2204.09596}, 2022.

\bibitem{oksendal2013stochastic}
B.~Oksendal, \emph{Stochastic differential equations: an introduction with
  applications}.\hskip 1em plus 0.5em minus 0.4em\relax Springer Science \&
  Business Media, 2013.

\bibitem{kappen2005path}
H.~J. Kappen, ``Path integrals and symmetry breaking for optimal control
  theory,'' \emph{Journal of statistical mechanics: theory and experiment},
  vol. 2005, no.~11, p. P11011, 2005.

\bibitem{williams2017model}
G.~Williams, A.~Aldrich, and E.~A. Theodorou, ``Model predictive path integral
  control: From theory to parallel computation,'' \emph{Journal of Guidance,
  Control, and Dynamics}, vol.~40, no.~2, pp. 344--357, 2017.

\bibitem{vrushabh2020robust}
D.~Vrushabh, P.~Akshay, K.~Sonam, S.~Wagh, and N.~M. Singh, ``Robust path
  integral control on stochastic differential games,'' in \emph{2020 28th
  Mediterranean Conference on Control and Automation (MED)}.\hskip 1em plus
  0.5em minus 0.4em\relax IEEE, 2020, pp. 665--670.

\bibitem{bacsar2008h}
T.~Ba{\c{s}}ar and P.~Bernhard, \emph{H-infinity optimal control and related
  minimax design problems: a dynamic game approach}.\hskip 1em plus 0.5em minus
  0.4em\relax Springer Science \& Business Media, 2008.

\bibitem{satoh2016iterative}
S.~Satoh, H.~J. Kappen, and M.~Saeki, ``An iterative method for nonlinear
  stochastic optimal control based on path integrals,'' \emph{IEEE Transactions
  on Automatic Control}, vol.~62, no.~1, pp. 262--276, 2016.

\bibitem{friedman1975stochastic}
A.~Friedman, \emph{Stochastic differential equations and applications, vol.
  1.}\hskip 1em plus 0.5em minus 0.4em\relax Academic Press, 1975.

\bibitem{patil2022chance}
A.~Patil, A.~Duarte, A.~Smith, F.~Bisetti, and T.~Tanaka, ``Chance-constrained
  stochastic optimal control via path integral and finite difference methods,''
  in \emph{2022 IEEE 61st Conference on Decision and Control (CDC)}.\hskip 1em
  plus 0.5em minus 0.4em\relax IEEE, 2022, pp. 3598--3604.

\bibitem{theodorou2010generalized}
E.~Theodorou, J.~Buchli, and S.~Schaal, ``A generalized path integral control
  approach to reinforcement learning,'' \emph{The Journal of Machine Learning
  Research}, vol.~11, pp. 3137--3181, 2010.

\bibitem{durrett2019probability}
R.~Durrett, \emph{Probability: theory and examples}.\hskip 1em plus 0.5em minus
  0.4em\relax Cambridge university press, 2019, vol.~49.

\bibitem{isaacs1999differential}
R.~Isaacs, \emph{Differential games: a mathematical theory with applications to
  warfare and pursuit, control and optimization}.\hskip 1em plus 0.5em minus
  0.4em\relax Courier Corporation, 1999.

\end{thebibliography}

\end{document}